\documentclass{article}
\begin{document}
\centerline{\large\bf On the Solution of Gauss Circle Problem Conjecture (Revised)}

\[
\]

\begin{quote}

\centerline{\bf Nikolaos D. Bagis \rm}
\centerline{Aristotle University of Thessaloniki}
\centerline{Thessaloniki, Greece}
\centerline{email: nikosbagis@hotmail.gr}
\[
\]

\centerline{\bf Abstract}
We give an asymptotic formula for the mean value of the number of representations of an integer as sum of two squares known as the Gauss circle problem.        
\[
\]
\textbf{keywords}: \textrm{Gauss Circle Problem; Conjecture; Proof; Sums of Squares; Asymptotic;}

\end{quote}

\section{Introduction}

In geometry the number of integer points lying in the circle
\begin{equation}
x^2+y^2=n,
\end{equation} 
is 'say' $r_2(n)$. The equivalent of the above proposition in number theory is saying, the number of representations of the integer $n$ as a sum of two squares is $r_2(n)$. Hence evaluating the function $r_2(n)$ has both geometrical and number theoretic meaning. Gauss itself prove that
\begin{equation}
\lim_{x\rightarrow\infty}\frac{1}{x}\sum_{n\leq x}r_2(n)=\pi .
\end{equation}

The problem we are describing, firstly began by finding the best possible constant $\theta>0$ of the asymptotic expansion
\begin{equation}
\sum_{n\leq x}r_2(n)=\pi x+O\left(x^{\theta+\epsilon}\right)\textrm{, }\forall\epsilon>0\textrm{, }x\rightarrow\infty.
\end{equation}
The first result it was that of Gauss which gave a first value $\theta=\frac{1}{2}$. For an analytic geometry proof of this case one can see [2].\\
After Gauss, many other scientists try to give better estimates with the most recent that of Huxley  (see [11]) in 2003:
$$
\sum_{n\leq x}r_2(n)=\pi x+O\left(x^{131/416}\right)\textrm{, }x\rightarrow\infty.
$$
The Gauss circle problem (conjecture) states that $\theta=\frac{1}{4}$ in (3).\\
There is a conversion procedure due to Richert that expresses the Gauss circle problem in terms of divisor problem of Dirichlet i.e
\begin{equation}
\Delta(x)=\sum_{n\leq x}d(n)-(\log x+2\gamma-1)x\textrm{, where }d(n)=\sum_{d|n}1.
\end{equation}    
The conjecture for the Dirchlet divisor problem is
\begin{equation}
\Delta(x)=O\left(x^{1/4+\epsilon}\right)\textrm{, }\forall \epsilon>0\textrm{, }x\rightarrow\infty.
\end{equation}
A result of Hardy shows that
\begin{equation}
\textrm{lim sup}_{x\rightarrow\infty}\frac{\Delta(x)}{x^{1/4}}=\infty
\end{equation}
and the complementary result of Ingham and Landau
\begin{equation}
{\textrm{lim inf}}_{x\rightarrow\infty}\frac{\Delta(x)}{x^{1/4}}=-\infty.
\end{equation} 
From analytic number theory point of view, Jacobi (see [4,12]) give a simple evaluation of the function $r_2(n)$. He defined the theta function (see [4,8]): $\vartheta(q)=\sum_{n\in\bf Z\rm}q^{n^2}$, for $|q|<1$, then he write 
$$
\vartheta^2(q)=\sum_{(n,m)\in\bf Z\times Z\rm}q^{n^2+m^2}=\sum^{\infty}_{n=0}r_2(n)q^n.
$$     
But as proved also by Jacobi (see details in [4]) holds the following theorem\\
\\
\textbf{Theorem 1.} (Jacobi)\\
If $q=e^{-\pi\sqrt{r}}$, $r>0$, then
\begin{equation}
\vartheta_3(q)=\sum^{\infty}_{n=-\infty}q^{n^2}=\sqrt{\frac{2K}{\pi}}.
\end{equation}

Hence from the equality
$$
\frac{2K}{\pi}={}_{2}F_1\left(\frac{1}{2},\frac{1}{2};1;k_r^2\right)=1+4\sum^{\infty}_{m=1}\frac{q^m}{1+q^{2m}}=1+4\sum^{\infty}_{m=1}q^m\sum^{\infty}_{l=0}(-1)^lq^{2ml}=
$$
\begin{equation}
=1+4\sum^{\infty}_{m=1}\sum^{\infty}_{l=0}(-1)^lq^{(2l+1)m}
\end{equation}
(the function $k_r$ being used here is called elliptic singular modulus and is defined in terms of Weber's functions (see [4,8]) as 
$k^2_r=16q\prod^{\infty}_{n=1}\left(\frac{1+q^{2n}}{1+q^{2n-1}}\right)^8$) writing the double series (9) under one sum using divisors sums he obtained the next result\\
\\
\textbf{Theorem 2.} (Jacobi)\\
For $n=1,2,\ldots$ we have 
\begin{equation}
r_2(n)=4\sum_{d-odd,\textrm{ }d|n}(-1)^{\frac{d-1}{2}}
\end{equation}
and $r(0)=1$.\\

Here we are interested only on the Gauss circle problem point of view (relation (3)) and not the Dirichlet divisor problem. The Jacobi formula (of Theorem 2) as it stands is very comfortable to use it directly and get representations of a nonnegative integer $n$ as sums of two squares. Next we give a precise formula for the mean value of $r_2(n)$, which is due to Hardy and Voronoi near 1904, (see [9,7]) 
\begin{equation}
\sum'_{n\leq x}r_2(n)=\pi x+x^{1/2}\sum^{\infty}_{n=1}\frac{r_2(n)}{\sqrt{n}}J_1(2\pi \sqrt{nx}),
\end{equation}
where the prime on summation means that the last term is multiplied with $1/2$. Also $J_{\nu}(z)$ is the general Bessel function of the first kind and $\nu$-th order and given by
\begin{equation}
J_{\nu}(z):=\sum^{\infty}_{k=0}\frac{(-1)^k}{k!\Gamma(k+\nu+1)}\left(\frac{z}{2}\right)^{\nu+2k}\textrm{, }0\leq |z|<\infty\textrm{, }\nu\in \textbf{C}.
\end{equation}

\section{A more simple statement of Gauss circle problem}

If $a,b\in\textbf{R}$, then we define
\begin{equation}
M_s(a,b)=\sum^{\infty}_{k-odd\textrm{, }k=1}(-1)^{\frac{k+1}{2}}\frac{\cos(a+b\sqrt{k})}{k^s}
\end{equation}
\begin{equation}
N_s(a,b)=\sum^{\infty}_{k-odd\textrm{, }k=1}(-1)^{\frac{k+1}{2}}\frac{\sin(a+b\sqrt{k})}{k^s}
\end{equation}
and
\begin{equation}
P_{s}(a,b)=\sum^{\infty}_{n=1}\frac{M_{s}(a,b\sqrt{n})}{n^s}\textrm{, }Q_{s}(a,b)=\sum^{\infty}_{n=1}\frac{N_{s}(a,b\sqrt{n})}{n^s}.
\end{equation}
\\
\textbf{Theorem 3.}\\
We have
$$
R(x)=\sum_{n\leq x}r_2(n)-x\pi=\frac{x^{1/4}}{\pi}P_{3/4}\left(\frac{\pi}{4},2\pi\sqrt{x}\right)+\sum^{N}_{s=1}\frac{(-1)^sc_1(2s)P_{s+3/4}\left(\frac{\pi}{4},2\pi\sqrt{x}\right)}{2^{4s}\pi^{2s+1}x^{s-1/4}}-
$$
\begin{equation}
-\sum^{N}_{s=0}\frac{(-1)^sc_1(2s+1)Q_{s+5/4}\left(\frac{\pi}{4},2\pi\sqrt{x}\right)}{2^{4s+2}\pi^{2s+2}x^{s+1/4}}+O\left(c_1(2N)4^{-N}x^{-N-1/2}\right)
\end{equation}
where $c_1(m)=(-1)^m\frac{\left(-\frac{1}{2}\right)_m\left(\frac{3}{2}\right)_m}{m!}$.\\
\\
\textbf{Proof.}\\
From (11) and Theorem 2 we have
$$
\sqrt{x}\sum_{n=1}^{\infty}\frac{r_2(n)}{\sqrt{n}}J_1(2\pi\sqrt{nx})=\sqrt{x}\lim_{N_1\rightarrow\infty}\sum_{n=1}^{N_1}\frac{r_2(n)}{\sqrt{n}}J_1(2\pi\sqrt{nx})=
$$
$$
\sqrt{x}\lim_{N_1\rightarrow\infty}\sum^{N_1}_{n=1}\left(\sum_{d-odd, d|n}(-1)^{\frac{d-1}{2}}\right)\frac{1}{\sqrt{n}}J_1(2\pi\sqrt{nx})=
$$
$$
\lim_{N_1\rightarrow\infty}\sqrt{x}\sum_{n(2m-1)\leq N_1}\frac{(-1)^m}{\sqrt{n(2m-1)}}J_1\left(2\pi\sqrt{n(2m-1)x}\right)=
$$
\begin{equation}
\sqrt{x}\lim_{N_1\rightarrow\infty}\sum_{p-odd, n p\leq N_1}\frac{(-1)^{\frac{p-1}{2}}}{\sqrt{np}}J_1\left(2\pi\sqrt{npx}\right)
\end{equation}
The function $J_1(x)$ has the following asymptotic expansion as $x\rightarrow\infty$
$$
J_1(x)\approx\sqrt{\frac{2}{\pi x}}[\cos\left(x-\frac{3\pi}{4}\right)\sum^{\infty}_{n=0}\frac{(-1)^nc_1(2n)}{(2x)^{2n}}-
$$
\begin{equation}
-\sin\left(x-\frac{3\pi}{4}\right)\sum_{n=0}^{\infty}\frac{(-1)^nc_1(2n+1)}{(2x)^{2n+1}}].
\end{equation}
The error due to stopping the summation at any term is the order of  magnitude of that term multiplied by $1/x$. Hence, using (18) in (17) we get (16).\\
\\

Setting $N=1$ in (16), we get
$$
R(x)=-\lim_{N_1\rightarrow\infty}\sum_{\scriptsize
\begin{array}{cc} 
	np\leq N_1\\
	p-odd
\end{array}\normalsize}[\frac{105 (-1)^{\frac{p+1}{2}}\sin\left(2\pi  \sqrt{npx}+\frac{\pi}{4}\right)}{4096\pi^3 (np)^{9/4} x^{5/4}}-
$$
$$
-\frac{15(-1)^{\frac{p+1}{2}} \cos \left(2\pi  \sqrt{npx}+\frac{\pi}{4}\right)}{256\pi^2 (np)^{7/4} x^{3/4}}
+\frac{3(-1)^{\frac{p+1}{2}} \sin\left(2\pi  \sqrt{npx}+\frac{\pi}{4}\right)}{8\pi (np)^{5/4} \sqrt[4]{x}}
$$
$$
-\frac{2 (-1)^{\frac{p+1}{2}} \sqrt[4]{x} \cos \left(2\pi\sqrt{np x}+\frac{\pi }{4}\right)}{(np)^{3/4}}]+O\left(x^{-3/4}\right),
$$
where $p=2l+1$.\\
\\
\textbf{Theorem 4.}\\
The Gauss circle problem reduces finding the rate of convergence of 
\begin{equation}
R(x)=\frac{1}{x^{1/4}}\left(\sum_{n\leq x}r_{2}(n)-\pi x\right),
\end{equation}
which is equivalent to that of
$$
S_{\infty}(x)=\lim_{N_1\rightarrow\infty}\sum_{n(2l-1)\leq N_1}\frac{(-1)^{l-1}\cos \left(2\pi\sqrt{n(2l-1)x}+\frac{\pi}{4}\right)}{(n(2l-1))^{3/4}}=
$$
\begin{equation}
=\lim_{N_1\rightarrow\infty}\sum^{N_1}_{n=1}\frac{r_2(n)\cos\left(2\pi\sqrt{nx}+\frac{\pi}{4}\right)}{n^{3/4}}.
\end{equation} 
If one manage to show that
\begin{equation}
S_{\infty}(x)=\sum^{\infty}_{n=1}\frac{r_2(n)\cos\left(2\pi\sqrt{nx}+\frac{\pi}{4}\right)}{n^{3/4}}
\end{equation}
is convergent and
\begin{equation}
S_{\infty}(x)=O\left(x^{\epsilon}\right)\textrm{, }x\rightarrow\infty,
\end{equation}
for all $\epsilon>0$, then the problem is solved.\\
\\

Instead of $S_{\infty}(x)$ we will study 
\begin{equation}
P_M(x)=\sum^{M}_{n=1}r_2(n)\frac{e^{2\pi i \sqrt{nx}}}{n^{3/4}}.
\end{equation}
It holds (this is a known-beter result of Gauss estimate):
\begin{equation}
\sum'_{0\leq n\leq t}r_2(n)=\pi t+O\left(t^{1/3}\right)\textrm{, }t\rightarrow\infty.
\end{equation}
Setting the function $g(x)$ to be such that (view also relation (11)):
\begin{equation}
\sum'_{0\leq n\leq t}r_2(n)=\pi t+g(t),
\end{equation}
we can write
\begin{equation}
g(x)= x^{1/4}S_{\infty}(x)=Re\left(e^{i\pi/4}x^{1/4}P_{\infty}(x)\right).
\end{equation}
Note that if $A,B$ are reals, then 
$$
\left|Re\left(e^{i \pi/4}(A+iB)\right)\right|=\frac{\sqrt{2}}{2}|A-B|\leq \sqrt{A^2+B^2}.
$$
Also we can write equivalent from (23):
$$
P_M(x)=\left(\sum^{M}_{n=1}r_2(n)\right)\frac{e^{2\pi i \sqrt{Mx}}}{M^{3/4}}-\int^{M}_{1}\left(\sum_{0\leq n\leq t}r_2(n)\right)\frac{d}{dt}\left(\frac{e^{2\pi i\sqrt{tx}}}{t^{3/4}}\right)dt=
$$
$$
=\left(\pi M+g(M)\right)\frac{e^{2\pi i \sqrt{Mx}}}{M^{3/4}}-\int^{M}_{1}\left(\pi t+g(t)\right)\frac{d}{dt}\left(\frac{e^{2\pi i \sqrt{tx}}}{t^{3/4}}\right)dt.
$$
Now using  
\begin{equation}
g(x)=O\left(\sqrt{x}\right)\textrm{, }x\rightarrow+\infty,
\end{equation}
we have
$$
P_M(x)=\left(\pi M+g(M)\right)\frac{e^{2\pi i \sqrt{Mx}}}{M^{3/4}}-\int^{M}_{1}\left(\pi t+g(t)\right)\frac{d}{dt}\left(\frac{e^{2\pi i \sqrt{tx}}}{t^{3/4}}\right)dt=
$$
$$
=\pi \left(M^{1/4}e^{2\pi i \sqrt{M x}}-\int^{M}_{1}t\frac{d}{dt}\left(\frac{e^{2\pi i \sqrt{tx}}}{t^{3/4}}\right)dt\right)
+e^{2 \pi i \sqrt{M x}}g(M)M^{-3/4}-
$$
$$
-\int^{M}_{1}g(t)\frac{d}{dt}\left(\frac{e^{2\pi i \sqrt{tx}}}{t^{3/4}}\right)dt.
$$
And we can  write
$$
|P_{M}(x)|\leq \pi \left|M^{1/4}e^{2\pi i \sqrt{M x}}-\int^{M}_{1}t\frac{d}{dt}\left(\frac{e^{2\pi i \sqrt{tx}}}{t^{3/4}}\right)dt\right|+
$$
\begin{equation}
+\left|\int^{M}_{1}g(t)\frac{d}{dt}\left(\frac{e^{2\pi i \sqrt{tx}}}{t^{3/4}}\right)dt\right|.
\end{equation}
Now the first of the two integrals when $M\rightarrow\infty$ is
$$
\int^{M}_{1}t\frac{d}{dt}\left(\frac{e^{2\pi i\sqrt{tx}}}{t^{3/4}}\right)dt=
$$
$$
=-e^{2\pi i \sqrt{x}}+\sqrt[4]{M}e^{2i\pi\sqrt{xM}}+\frac{1-i}{\sqrt[4]{x}}\textrm{Erfi}\left((1+i)\sqrt[4]{x}\sqrt{\pi}\right)-
$$
$$
-\frac{1-i}{\sqrt[4]{x}}\textrm{Erfi}\left((1+i)\sqrt[4]{xM}\sqrt{\pi}\right).
$$
One can easily see that the first term in the right side of (28) is always bounded function of $x$, when we take the limit $M\rightarrow \infty$. The second term can be evaluated with integration by parts as follows 
$$
\left|\int^{M}_{1}g(t)\frac{d}{dt}\left(\frac{e^{2\pi i \sqrt{tx}}}{t^{3/4}}\right)dt\right|
=\left|\left[g(t)\frac{e^{2\pi i \sqrt{t x}}}{t^{3/4}}\right]^{M}_{1}-\int^{M}_{1}g'(t)\frac{e^{2\pi i \sqrt{tx}}}{t^{3/4}}dt\right|\leq
$$
$$
\leq C_1+\left|\int^{M}_{1}g'(t)\frac{e^{2\pi i \sqrt{tx}}}{t^{3/4}}dt\right|.\eqno{(28.1)}
$$
Hence (when $|g'(t)|\leq Ct^{1/4+\epsilon}$), in view of the Lemmas 1,2 below we have $\lim_{x\rightarrow\infty}\left|P_{\infty}(x)\right|\leq C_1$ when $\epsilon<0$ and in general if $\epsilon>0$, then
$$
\left|\int^{\infty}_{1}g'(t)\frac{e^{2\pi i \sqrt{tx}}}{t^{3/4}}dt\right|=\infty.
$$
Hence $\left|P_{\infty}(x)\right|=O(1)$, $x>>1$, for every $\epsilon<0$. But this is known to be false. Hence we have $\epsilon>0$. Set now
$$
I_f(M,x):=\int^{\infty}_{1}f'(t)\frac{e^{2\pi i \sqrt{tx}}}{t^{3/4}}dt.
$$
It holds the following:\\
\\
\textbf{Proposition 1.}\\ When $x>>1$, we have that
\begin{equation}
\lim_{M\rightarrow\infty}\left|\int^{M}_{1}\frac{g'(t)e^{2\pi i \sqrt{tx}}}{t^{3/4}}dt\right|,
\end{equation}
is not bounded function of $x$ and $|g'(x)|$ is not monotone.\\
\\
\textbf{Proof.}\\
If $I_{g}(\infty,x)$ was bounded then $\left|P_{\infty}(x)\right|$ will be bounded, which is not true (see [2]).\\ 
We will show that if $\left|g'(t)\right|$ is monotone then $P_M(x)$ is bounded always and this leads to contradiction.\\
\textbf{i)} If $\left|g'(t)\right|$ where increasing, then for an arbitrary $y>1$ will have that exists $\xi$ such that  $1<y\leq\xi\leq 2y$ and
$$
\left|g'(y)\right|\leq\left|g'(\xi)\right|=\left|\frac{g(2y)-g(y)}{y}\right|\leq\left|g'(2y)\right|.
$$
Hence form Gauss estimate we get
$$
\left|g'(y)\right|\leq\frac{|g(2y)|+|g(y)|}{y}\leq \frac{A}{\sqrt{y}}\textrm{, }y>>1, 
$$
which is contradiction from (28.1) and $\left|P_{\infty}(x)\right|\neq O(1)$, $x\rightarrow \infty$.\\
\textbf{ii)} If $\left|g'(t)\right|$ where decreasing, then for arbitrary $y>1$ will exist $\xi$ with $1<\xi<y$ such that
$$
\left|g'(y)\right|\leq\left|g'(\xi)\right|=\left|\frac{g(y)-g(1)}{y-1}\right|\leq |g'(1)|.
$$
But from Gauss estimate exists $B\geq 1$ such that $|g(y)|\leq B \sqrt{y}$ and 
$$
|g'(y)|\leq\frac{|g(y)-g(1)|}{|y-1|}\leq \frac{B'}{\sqrt{y}}\textrm{, when }y>>1.
$$
But this is contradiction as in proof of (i) above.\\Hence from (i),(ii), we get that if $\left|g'(x)\right|$ were monotone then (28.1) will be bounded. Hence $P_{\infty}(x)$ will be bounded, which is not true.\\
\\
\textbf{Proposition 2.}\\
It is true that
$$
\left|g(x)\right|=O\left(x^{1/4}\left|\int^{\infty}_{1}g'(t)\frac{e^{2\pi i \sqrt{tx}}}{t^{3/4}}dt\right|\right)\textrm{, }x\rightarrow+\infty. \eqno{(30)}
$$
\textbf{Proof.}\\
It follows from (26) and the above analysis.\\
\\
\textbf{Lemma 1.}\\
If $1>\delta>1/2$, then 
$$
\int^{\infty}_{1}\frac{e^{2\pi i \sqrt{tx}}}{t^{\delta}}dt\eqno{(30.1)}
$$
converges uniformly but not absolutely to a function of $x>0$. Also
\begin{equation}
\lim_{x\rightarrow\infty}\left|\int^{\infty}_{1}\frac{e^{2\pi i \sqrt{tx}}}{t^{\delta}}dt\right|=0.
\end{equation}
Moreover for every $1>\delta>1/2$ and for every $C\in\textbf{R}$, we have that the improper integral
\begin{equation}
\int^{\infty}_{1}(\log t)^C\frac{e^{2\pi i\sqrt{tx}}}{t^{\delta}}dt
\end{equation}
converges uniformly but not absolutely to a function of $x>0$ and
\begin{equation}
\lim_{x\rightarrow+\infty}\left|\int^{\infty}_{1}(\log t)^C\frac{e^{2\pi i\sqrt{tx}}}{t^{\delta}}dt\right|=0.
\end{equation}
\\
\textbf{Proof.}\\
Assume that $1>\delta>1/2$, then
$$
J_1(M,\delta,x):=\int^{M}_{1}\frac{e^{2\pi i\sqrt{tx}}}{t^{\delta}}dt=\int^{\sqrt{M}}_{1}\frac{e^{2\pi i t\sqrt{x}}}{t^{2\delta}}2tdt=
$$
$$
=2\int^{\sqrt{M}}_{1}\frac{e^{2 \pi i t\sqrt{x}}}{t^{2\delta-1}}dt.
$$
Since $\delta>1/2$ we set $\epsilon=2\delta-1>0$ and get easily 
$$
\int^{\infty}_{1}\frac{e^{2\pi i t\sqrt{x}}}{t^{\epsilon}}dt=\int^{\infty}_{1}\frac{\cos\left(2\pi  t\sqrt{x}\right)}{t^{\epsilon}}dt+i\int^{\infty}_{1}\frac{\sin\left(2\pi t\sqrt{x}\right)}{t^{\epsilon}}dt=
$$
\begin{equation}
=E(\epsilon,-2 i \pi \sqrt{x})\textrm{, }\epsilon, x>0,
\end{equation}
where 
\begin{equation}
E(\nu,z)=\int^{\infty}_{1}\frac{e^{-zt}}{t^{\nu}}dt,
\end{equation}
is the exponential integral $E$. This function is defined for all $z\in\textbf{C}-\textbf{R}_{<0}$. Hence in our case 
$J_1(\infty,\delta,x)=2E(2\delta-1,-2i\pi\sqrt{x})<\infty$. For to evaluate $\lim_{x\rightarrow+\infty}|J_1(\infty,\delta,x)|$, we integrate by parts
$$
|J_1(M,\delta,x)|=\left|2\int^{\sqrt{M}}_{1}\frac{d}{dt}\left(\frac{e^{2\pi i t \sqrt{x}}}{2\pi i \sqrt{x}}\right)\frac{dt}{t^{2\delta-1}}\right|=
$$
$$
=\left|\frac{1}{\pi i\sqrt{x}}\left[\frac{e^{2\pi i t\sqrt{x}}}{t^{2\delta-1}}\right]^{t=\sqrt{M}}_{t=1}-\frac{1}{\pi i\sqrt{x}}(1-2\delta)\int^{\sqrt{M}}_{1}\frac{e^{2\pi i t \sqrt{x}}}{t^{2\delta}}dt\right|\leq
$$
$$
\leq\left|\frac{1}{\pi \sqrt{x}}\left[\frac{e^{2\pi i t\sqrt{x}}}{t^{2\delta-1}}\right]^{t=\sqrt{M}}_{t=1}\right|+\frac{1}{\pi \sqrt{x}}(2\delta-1)\left|\int^{\sqrt{M}}_{1}\frac{e^{2\pi i t \sqrt{x}}}{t^{2\delta}}dt\right|$$
Now since $1>\delta>1/2$ we can take the limits first  $M\rightarrow\infty$ and then $x\rightarrow\infty$ to get the result.\\
For the second integral (31) assume that $C\geq 1$ (the case $C<1$ is easy). Setting $\epsilon_{1,2}>0$ to be such that $\delta=\frac{1}{2}+\epsilon_1+\epsilon_2$, we get
$$
J_2(C,M,\delta,x)=\int^{M}_{1}(\log t)^C\frac{e^{2\pi i \sqrt{tx}}}{t^{\delta}}dt=\int^{M}_{1}\frac{(\log t)^C}{t^{\epsilon_1}}\frac{e^{2\pi i \sqrt{tx}}}{t^{1/2+\epsilon_2}}dt.
$$
Now the function $(\log t)^C/t^{\epsilon_1}$ is decreasing to $0$, for all $t>t_0>0$, where $t_0$ depending on $C,\epsilon_1$. Also from the first question the integral
$$
\int^{M}_{1}\frac{e^{2\pi i\sqrt{tx}}}{t^{1/2+\epsilon_2}}dt=2\int^{\sqrt{M}}_{1}\frac{e^{2\pi i t\sqrt{x}}}{t^{2\epsilon_2}}dt<\infty
$$ 
is bounded when $M\rightarrow\infty$. Hence from Dirichlet's test for integrals, we have that $J_2(C,\infty,\delta,x)$ converges uniformly. Now write
$$
J_2(C,M,\delta,x)=2^{C+1}\int^{\sqrt{M}}_{1}\frac{(\log t)^C}{t^{2\epsilon_1}}\frac{e^{2\pi i t\sqrt{x}}}{t^{2\epsilon_2}}dt=
$$
$$
=2^{C+1}\int^{\sqrt{M}}_{1}\frac{(\log t)^{C}}{t^{2\epsilon_1+2\epsilon_2}}\frac{d}{dt}\left(\frac{e^{2\pi i t\sqrt{x}}}{2\pi i\sqrt{x}}\right)dt=
$$
\begin{equation}
=\frac{2^{C+1}}{2\pi i\sqrt{x}}\left[\frac{(\log t)^{C}}{t^{2\epsilon_1+2\epsilon_2}}e^{2\pi i t\sqrt{x}}\right]^{t=\sqrt{M}}_{t=1}-\frac{2^{C+1}}{2\pi i\sqrt{x}}\int^{\sqrt{M}}_{1}\frac{d}{dt}\left(\frac{(\log t)^{C}}{t^{2\epsilon_1+2\epsilon_2}}\right)e^{2\pi i t\sqrt{x}}dt.
\end{equation}
But
\begin{equation}
\int^{\infty}_{1}\frac{d}{dt}\left(\frac{(\log t)^{C}}{t^{2\epsilon_1+2\epsilon_2}}\right)e^{2\pi i t\sqrt{x}}dt=\int^{\infty}_{1}\frac{(\log t)^{C-1}(C-2(\epsilon_1+\epsilon_2)\log t)}{t^{2\epsilon_1+2\epsilon_2+1}}e^{2\pi i t\sqrt{x}}dt
\end{equation}
and the function
$$
\frac{(\log t)^{C-1}(C-2(\epsilon_1+\epsilon_2)\log t)}{t^{2\epsilon_1+2\epsilon_2+1}},
$$
belongs to $L^{1}[1,+\infty)$. Hence from (35) and (36) and applying Riemann-Lebesgue theorem we get $\lim_{x\rightarrow\infty}J_2(C,\infty,\delta,x)=0$ and hence we get (32).\\
\\
\textbf{Remarks.}\\
If $\delta$ is any positive real number and $M>1$, $x>0$, we set (as above)
\begin{equation}
J_1(M,\delta,x)=\int^{M}_{1}\frac{e^{2\pi i \sqrt{tx}}}{t^{\delta}}dt
\end{equation}
and more general
\begin{equation}
J_2(C,M,\delta,x)=\int^{M}_{1}(\log(t))^C\frac{e^{2\pi i \sqrt{tx}}}{t^{\delta}}dt.
\end{equation}
Then
\begin{equation}
J_1(M,\delta,x)=\frac{1}{\pi i\sqrt{x}}\left[\frac{e^{2\pi i t\sqrt{x}}}{t^{2\delta-1}}\right]^{t=\sqrt{M}}_{t=1}-\frac{1}{\pi i\sqrt{x}}\left(\frac{1}{2}-\delta\right)\cdot  J_1\left(M,\frac{1}{2}+\delta,x\right)
\end{equation}
and
$$
J_2(C,M,\delta,x)=\frac{2^{C+1}}{2\pi i\sqrt{x}}\left[(\log t)^C\frac{e^{2\pi i t\sqrt{x}}}{t^{2\delta-1}}\right]^{t=\sqrt{M}}_{t=1}-
$$
\begin{equation}
-\frac{1}{\pi i\sqrt{x}}CJ_2\left(C-1,M,\frac{1}{2}+\delta,x\right)+\frac{1}{2\pi i\sqrt{x}}(2\delta-1)J_2\left(C,M,\frac{1}{2}+\delta,x\right).
\end{equation}
Moreover for every $\delta>1/2$ and for every fixed $C\in\textbf{R}$, we have that
\begin{equation}
\int^{\infty}_{1}(\log t)^C\frac{e^{2\pi i\sqrt{tx}}}{t^{\delta}}dt,
\end{equation}
exists for every $x>0$ and
\begin{equation}
\lim_{x\rightarrow+\infty}\int^{\infty}_{1}(\log t)^C\frac{e^{2\pi i\sqrt{tx}}}{t^{\delta}}dt=0.
\end{equation}
\\
\textbf{Lemma 2.}\\
If $C\in\textbf{R}$ and assume that 
$$
h(t)=t^{1/4+\epsilon}(\log t)^C,
$$ 
then\\ 
\textbf{1)} If $\epsilon<0$  we have
$$
\lim_{x\rightarrow+\infty}\left|\int^{\infty}_{1}h(t)\frac{e^{2\pi i \sqrt{t x}}}{t^{3/4}}dt\right|=0
$$
and\\ 
\textbf{2)} If $\epsilon>0$ we have 
$$
\lim_{x\rightarrow+\infty}\left|\int^{\infty}_{1}h(t)\frac{e^{2\pi i \sqrt{t x}}}{t^{3/4}}dt\right|=+\infty.
$$
\\
\textbf{Proof.}\\
Easy. Use Lemma 1 and its remarks.\\
\\
\textbf{Theorem 5.}\\
If we assume that $g'(t)$ have asymptotic development of the form
$$
g'\left(e^t\right)\approx e^{t/4+\epsilon t}\sum^{\infty}_{n=0}A_nt^{C_n}\textrm{, }t\rightarrow+\infty.\eqno{(42.1)}
$$
Then
\begin{equation}
g'(t)\approx t^{1/4+\epsilon}\sum^{\infty}_{n=0}A_n\exp\left(C_n\log\log t\right)\textrm{, }t\rightarrow +\infty
\end{equation} 
and\\
\textbf{1)} If $\epsilon<0$ we have
$$
\lim_{x\rightarrow+\infty}\left|\int^{\infty}_{1}g'(t)\frac{e^{2\pi i \sqrt{t x}}}{t^{3/4}}dt\right|=0.
$$
\textbf{2)} If $\epsilon> 0$ we have
$$
\lim_{x\rightarrow+\infty}\left|\int^{\infty}_{1}g'(t)\frac{e^{2\pi i \sqrt{t x}}}{t^{3/4}}dt\right|=+\infty.
$$
\\
\textbf{Proof.}\\
Assume $g'(t)$ is of the form (43), then 
$$
g'\left(t\right)\frac{e^{2\pi i  \sqrt{tx}}}{t^{3/4}}=\sum^{\infty}_{n=0}A_nt^{1/4+\epsilon}(\log t)^{C_n}\frac{e^{2\pi i \sqrt{tx}}}{t^{3/4}}.
$$
Hence integrating we get
$$
\int^{M}_{1}g'\left(t\right) \frac{e^{2\pi i \sqrt{tx}}}{t^{3/4}}dt=\sum^{\infty}_{n=0}A_n\int^{M}_{1}t^{1/4+\epsilon}(\log t)^{C_n}\frac{e^{2\pi i \sqrt{tx}}}{t^{3/4}}dt.
$$
Hence when $\epsilon<0$, we get
$$
\lim_{x\rightarrow+\infty}\int^{\infty}_{1}g'(t)\frac{e^{2\pi i\sqrt{tx}}}{t^{3/4}}dt=
$$
$$
=\sum^{\infty}_{n=0}A_n\lim_{x\rightarrow+\infty}\int^{\infty}_{1}t^{1/4+\epsilon}(\log t)^{C_n}\frac{e^{2\pi i\sqrt{tx}}}{t^{3/4}}dt=\ldots=0
$$
and when $\epsilon>0$ we have
$$
\lim_{x\rightarrow+\infty}\left|\int^{\infty}_{1}g'(t)\frac{e^{2\pi i \sqrt{tx}}}{t^{3/4}}dt\right|=+\infty.
$$
By this way the theorem is proved.\\
\\
\textbf{Theorem 6.}\\
The function $g(x)$ have an asymptotic development of the form 
\begin{equation}
g(x)\approx -\sum^{\infty}_{n=0}A_n\Gamma\left(C_n+1,a\log x\right)a^{-C_n-1}\textrm{, }C_n=-n^{3/4},
\end{equation}
where $a=-5/4-\epsilon$, $\epsilon>0$ and $\Gamma(s,z)=\int^{\infty}_{z}t^{s-1}e^{-t}dt$.\\
Then if this asymptotic series is differentiatable, we have  
\begin{equation}
g'(x)\approx x^{1/4+\epsilon}\sum^{\infty}_{n=0}A_n\exp\left(-n^{3/4}\log \log x\right)\textrm{, }x\rightarrow+\infty
\end{equation}
and the conjecture is true.\\
\\
\textbf{Proof.}\\
From the asymptotic formula
\begin{equation}
-\Gamma\left(1-n^{3/4},a\log x\right)\left(\frac{1}{a\log x}\right)^{1-n^{3/4}}=c(a,x)+o(1)\textrm{, }n\rightarrow +\infty
\end{equation}
and (18), we get
\begin{equation}
x\log(\log x) \frac{x^{1/2}\frac{J_1\left(2\pi\sqrt{n x}\right)}{\sqrt{n}}}{-\Gamma(C_n+1,a\log x)\left(\frac{1}{a\log x}\right)^{C_n+1}}=c_1(a,x)+o(1)\textrm{, }n\rightarrow\infty,
\end{equation}
where $C_n=-n^{3/4}$. Hence 
$$
\frac{x^{1/2}\frac{J_1\left(2\pi\sqrt{n x}\right)}{\sqrt{n}}}{-\Gamma(C_n+1,a\log x)a^{-C_n-1}}=b(n,a,x)\frac{(\log x)^{|C_n|-1}}{x \log \log x},
$$
where $b(n,a,x)$ is a sequence such that 
$$
b(n,a,x)=c_1(a,x)+h_{n}(a,x)
$$ 
and
$$
\lim_{n\rightarrow+\infty}h_{n}(a,x)=0.
$$
Hence we can write
$$
g(x)=-\sum^{\infty}_{n=0}r_2(n)\Gamma(C_n+1,a\log x)a^{-C_n-1}\left(\frac{b(n,a,x)}{x\left(\log x\right)^{C_n+1}\log \log x}\right).
$$
Also for $x\geq e^{R_n}=e^{n^{3/4}-1}$, we have
$$
\frac{e^{R_n\log\log x}}{x\log \log x}=\frac{e^{R_n\log y-y}}{\log y}\leq \frac{e^{R_n\log R_n-R_n}}{R_n}\textrm{, }R_n=n^{3/4}-1.
$$
Hence we can write
$$
g(x)\approx  -\sum^{\infty}_{n=1}r_2(n)\Gamma\left(C_n+1,a\log x\right)a^{-C_n-1}\frac{c_1(a,x)e^{R_n\log \log x}}{x\log \log x}
$$
and
$$
g(x)\approx -\sum^{\infty}_{n=0}A_n\Gamma\left(C_n+1,a\log x\right)a^{-C_n-1},
$$
with 
$$
A_n=r_2(n)\frac{e^{R_n\log R_n-R_n}}{R_n}.
$$
Moreover if this asymptotic expansion is differentiatable, then, $\epsilon>0$ and
$$
g'(x)\approx x^{1/4+\epsilon}\sum^{\infty}_{n=0}A_n\exp\left(-n^{3/4}\log \log x\right)\textrm{, }x\rightarrow \infty
$$
and from Proposition 1 and Theorem 5 we get the conjecture.\\
\\
\textbf{Notes.}\\
\textbf{i)} The symbol litle$-o$ means that if $a_n=b_n+o(1)$, then exists positive null sequence $\epsilon_n$ such that $a_n=b_n+\epsilon_n$.\\
\textbf{ii)} As someone can see there is some normality using the estimate $g'(t)=O\left(t^{1/4+\epsilon}\right)$, $t\rightarrow\infty$ and the Lemmas 1,2 and Theorem 6. However when one uses the first derivative of (11) (assuming that the series is asymptoticaly differentiatable) we get that
$$
g'(x)\approx \pi\sum^{\infty}_{n=1}r_2(n)J_{0}\left(2\pi \sqrt{nx}\right)\textrm{, }x\rightarrow\infty ,
$$
which under using the asymptotic formula
$$
J_0(z)\approx \sqrt{\frac{2}{\pi z}}\left[\cos\left(z-\frac{\pi}{4}\right)+e^{\left|Im(z)\right|}O\left(|z|^{-1}\right)\right]\textrm{, }|z|\rightarrow\infty,
$$ 
results ($\epsilon>0$):
$$
g'(t)=O\left(t^{-1/4-\epsilon}\right)\textrm{, }t\rightarrow\infty,
$$
which also may occur. However in this case Lemmas 1,2 not hold when $\epsilon$ changes sign.
\[
\]

\centerline{\bf References}\vskip .2in

\noindent

[1]: M. Abramowitz and I.A. Stegun. 'Handbook of Mathematical Functions'. Dover Publications, New York. (1972)

[2]: G.E. Andrews, Number Theory. Dover Publications, New York, (1994)  

[3]: T. Apostol. 'Introduction to Analytic Number Theory'. Springer Verlag, New York, Berlin, Heidelberg, Tokyo, (1974)

[4]: J.V. Armitage W.F. Eberlein. 'Elliptic Functions'. Cambridge University Press. (2006)

[5]: N. Bagis. 'Some New Results on Sums of Primes'. Mathematical Notes, (2011), Vol. 90, No. 1, pp 10-19.

[6]: N. Bagis. 'Some Results on Infinite Series and Divisor Sums'.\\ arXiv:0912.48152v2 [math.GM] (2014).

[7]: George E. Andrews and Bruce C. Berndt. 'Ramanujan`s Lost Notebook IV'. Springer, New York, Heilderberg, Dordrecht, London. (2013)

[8]: J.M. Borwein and P.B. Borwein. 'Pi and the AGM'. John Wiley and Sons, Inc. New York, Chichester, Brisbane, Toronto, Singapore. (1987)

[9]: G.H. Hardy. 'Ramanujan Twelve Lectures on Subjects Suggested by his Life and Work, 3rd ed.' Chelsea. New York, (1999)

[10]: M.N. Huxley. 'Exponential Sums and Lattice Points III'. Proc. London Math. Soc. 87, 5910-609, (2003)  

[11]: C.G.J. Jacobi. 'Fundamenta Nova Functionum Ellipticarum'. Werke I, 49-239. (1829) 

[12]: Konrad Knopp. 'Theory and Applications of Infinite Series'. Dover Publications, Inc. New York. (1990).

[13]: William J. LeVeque. 'Fundamentals of Number Theory'. Dover Publications. New York. (1996)

[14]: E.T. Whittaker and G.N. Watson. 'A course on Modern Analysis'. Cambridge U.P. (1927)

[15]: H. Iwaniec and E. Kowalski. 'Analytic Number Theory'. American Mathematical Society. Colloquium Publications. Vol 53. Providence, Rhode Island. (2004)

[16]: Bruce C. Berndt, Sun Kim, Alexandru Zaharescu.'Circle and divisor problems, and double series of Bessel functions'. Adv.Math. 236. 24-59 (2013).

[17]: G.H. Hardy, E.M. Wright.'An Introduction to the Theory of Numbers'. Sixth edition. Oxford University Press. Oxford, (2008). 

[18]: N.G. de Bruijn. 'Asymptotic Methods in Analysis'. Dover Publications, Inc. New York. (1981).

[19]: E.C. Titchmarsh. 'The Theory of The Riemann Zeta-Function'. Oxford. At the Clarendon Press. (1951).

\end{document}